\documentclass[letterpaper, 10 pt, conference]{ieeeconf}  

\IEEEoverridecommandlockouts                              

\pdfoutput=1

\usepackage{graphicx}

\title{\LARGE \bf
Experimental Determination of Demand Response Control Models and Cost of Control for Ensembles of Window-Mount Air Conditioners}

\author{Drew A. Geller$^{1}$ and Scott Backhaus$^{2}$
\thanks{This work was supported by the Microgrid Program of the Office of Electricity within the U.S. Department of Energy}
\thanks{$^{1}$Drew Geller is with the Gas Transfer Systems Group, Los Alamos National Laboratory,
	Los Alamos, NM 87544, USA
        {\tt\small dgeller@lanl.gov}}%
\thanks{$^{2}$Scott Backhaus is with the Information Systems and Modeling Group, Los Alamos National Laboratory,
        Los Alamos, NM 87544, USA
        {\tt\small backhaus@lanl.gov}}%
}

\begin{document}

\maketitle
\thispagestyle{empty}
\pagestyle{empty}

\begin{abstract}
Control of consumer electrical devices for providing electrical grid services is expanding in both the scope and the diversity of loads that are engaged in control, but there are few experimentally-based models of these devices suitable for control designs and for assessing the cost of control. A laboratory-scale test system is developed to experimentally evaluate the use of a simple window-mount air conditioner for electrical grid regulation services.  The experimental test bed is a single, isolated air conditioner embedded in a test system that both emulates the thermodynamics of an air conditioned room and also isolates the air conditioner from the real-world external environmental and human variables that perturb the careful measurements required to capture a model that fully characterizes both the control response functions and the cost of control. The control response functions and cost of control are measured using harmonic perturbation of the temperature set point and a test protocol that further isolates the air conditioner from low frequency environmental variability.
\end{abstract}

\section{INTRODUCTION}
Renewable sources of electrical power such as solar and wind farms may have large output variations over a wide range of time scales, from minutes to hours to days. In the case of solar, the temporal variation in solar energy available is not only due to the daily cycle of the sun but also due to cloud cover and aerosols in the atmosphere.  The variability from additional renewables combines with the existing variability of loads to increase the need for power regulating services for continuously balancing generation and load over this wide range of time scales. Bulk electrical system (BES) operators already purchase or otherwise acquire grid ancillary services to provide this real power regulation. Examples include primary frequency regulation that operates on time scales from subsecond to ten seconds; secondary frequency regulation operating on time scales of a few seconds to ten minutes and spinning reserves or tertiary frequency regulation operating on time scales of ten minutes to an hour \cite{kundur1994power}.

Although the time scales of these services are quite wide, they are for the most part provided by traditional, central station synchronous generators. However, as the penetration of renewable resources increases, the energy supply from traditional generators is needed less and they are operated less making them unavailable to provide ancillary services for controlling the electrical grid. Advanced control of the renewable generators may be used to restore some of the controllable capacity \cite{Roald_2015} \cite{Roald_2016}.  However, new sources of controllable capacity will likely be needed at higher levels of renewable resource penetration \cite{Ela_2011}.

One concept for providing such regulating service involves controlling an ensemble of thermostatically controlled loads (TCLs) such as air conditioners (AC).  Because typical buildings or offices have relatively large heat capacities compared to the air contained in these spaces, it should be possible to make small changes in thermostat settings without affecting the comfort of the inhabitants of these spaces, at least over short periods of time.  By distributing these set point changes over a large number of TCLs, a large amount of power consumption can be controlled for the purpose of providing grid ancillary services, without strongly impacting any particular load for a sustained period.

Although TCLs, e.g. large commercial heating, ventilation and air conditioning (HVAC) systems and ensembles of many small residential AC systems, are attractive control targets for providing ancillary services, there are few experimentally-based models of these devices suitable for control designs and for assessing the cost of control. A review of both experimental and simulation-based evaluations of control models and control for large commercial HVAC systems is given in \cite{Beil_2016}.  Here, we provide a discussion of the references in \cite{Beil_2016} that are most relevant to current work and discuss a few more recent references that are not included in \cite{Beil_2016}.

Experimental system identification was used in \cite{goddard} to determine an open-loop demand response control model that accounted for the nonlinear response of a large commercial HVAC system. This model was developed for the entire HVAC system in a 300,000 ft$^2$ building for time scales from 15-30 minutes. Its performance was tested \cite{Beil_2016} against historical PJM RegA signals \cite{PJM1}, which have a typical time scale of 5 minutes. The control performance was found to be adequate. This same open-loop control and commercial HVAC system combination was also used to experimentally determine the cost of control.  Transient testing \cite{beil} was used to determine excess energy consumption to return the HVAC system to its nominal internal state following a control action. In \cite{Beil_2016}, the open-loop control and HVAC system tracked a historical PJM RegA signal \cite{PJM1} for many hours and the average energy consumption was compared to days when no control tracking was performed.  In both cases, measurements showed significant excess energy consumption compared to the power regulation ancillary service provided.

Related experimental work was performed in \cite{lin} where frequency sweep methods were used to extract a plant transfer function for a single ventilation and supply air fan that was embedded in a larger commercial HVAC system.  The plant transfer function was used to design a closed-loop control for tracking both PJM RegA and RegD \cite{PJM1} frequency regulation signals. In both cases, the experimentally-driven design provided excellent tracking.  However, the authors of \cite{lin} did not attempt to measure the cost of control.

The work in \cite{goddard}\cite{beil}\cite{Beil_2016} and in \cite{lin} was performed on large commercial HVAC systems of similar configuration. They are comprised of a central chiller plant and chilled water loop that supplies air-to-water heat exchangers in a few air handling units (AHU). Large fans in the AHUs supply air to variable air volume controllers that ultimately supply cooling air to the individual conditioned spaces. This HVAC configuration is one of several typical layouts used in large commercial HVAC systems.

Related experimental work has also been performed on residential-scale HVAC systems \cite{Vrettos_2016}\cite{Vrettos_2016_2}, although the configuration of the tested HVAC systems may not be representative of those found in the existing residential building stock. In \cite{Vrettos_2016}, independent models of the thermal behavior of the residential building and the ventilation fans are formulated and fit to experimental data. These models are then incorporated into a model predictive control (MPC) framework that enables both long time-scale energy and frequency regulation capacity scheduling and short time-scale control to provide frequency regulation control. Experimental testing of the controls is performed in \cite{Vrettos_2016_2}, and the tracking of fast PJM RegD signals \cite{PJM1} was found to be excellent.

The authors of \cite{Vrettos_2016_2} did try to measure the cost of control while tracking PJM RegD signals. In contrast to \cite{beil} and \cite{Beil_2016}, they did not find any measurable effect. However, there are several potentially important differences between these experiments. First, the RegD reference signal used in \cite{Vrettos_2016_2} has a much faster time scale ($\approx$ 10 seconds) as compared to the RegA reference signal used in \cite{Beil_2016}. If the extra losses induced by frequency regulation control of the HVAC system are time scale dependent, they may be more sensitive to the slower variations in the RegA reference signal used in \cite{Beil_2016}. Second, the configuration of the HVAC system in \cite{Vrettos_2016_2} is not typical of residential HVAC systems and includes a large thermal mass in the cooling water loop (which is also not typical of residential systems).  These differences make it difficult to translate the results to HVAC systems in the existing building stock. Finally, the authors of \cite{Vrettos_2016_2} executed a very limited number of experiments where they attempted to measure the cost of control. Although the experimental system in \cite{Vrettos_2016_2} includes side-by-side HVAC systems to simultaneously measure the energy consumption with and without control, our recent experimental experience (described later in this manuscript) suggest that careful experimental methods and procedures are required to accurately determine the cost of control.

The work described in this manuscript attempts to remedy many of the issues discussed above. First, we utilize a very typical, low-cost window-mount AC purchased from a local retailer. Other than replacing the standard bulb-type thermostat with a remotely controlled relay and flexible logic programmed into the data acquisition system, the window-mount AC is used as purchased. Second, the window-mount AC is embedded in a test system that both emulates the thermodynamics of an air conditioned room and also isolates the AC from the real-world external environmental and human variables that perturb the careful measurements required to capture a model that fully characterizes both control response functions and the cost of control. Third, the control response functions and cost of control are measured using harmonic perturbation of the temperature set point and are explored over a wide range of perturbation time scales. Finally, the test protocol is cycled on and off to continually extract power consumption baselines from the same window-mount AC that is under test, which further isolates the AC from low frequency environmental variability.

The remainder of this manuscript is organized as follows. Section~\ref{sec:testbed} describes the experimental test bed and the necessary experimental precautions and procedures required to make accurate measurements of the control response functions and cost of control.  Section~\ref{sec:results} describes the measurements and the analysis methods and presents the main experimental results. Section~\ref{sec:discussion} discusses the results, places them in the context of previous work, and describes the implications of using window-mount ACs for demand response frequency regulation. Finally, Section~\ref{sec:conclusions} draws conclusions and discusses a path forward for future work. 
\section{Experimental Test Bed}\label{sec:testbed}
\begin{figure*}[!t]
	\centering
	\includegraphics[width=0.8\textwidth]{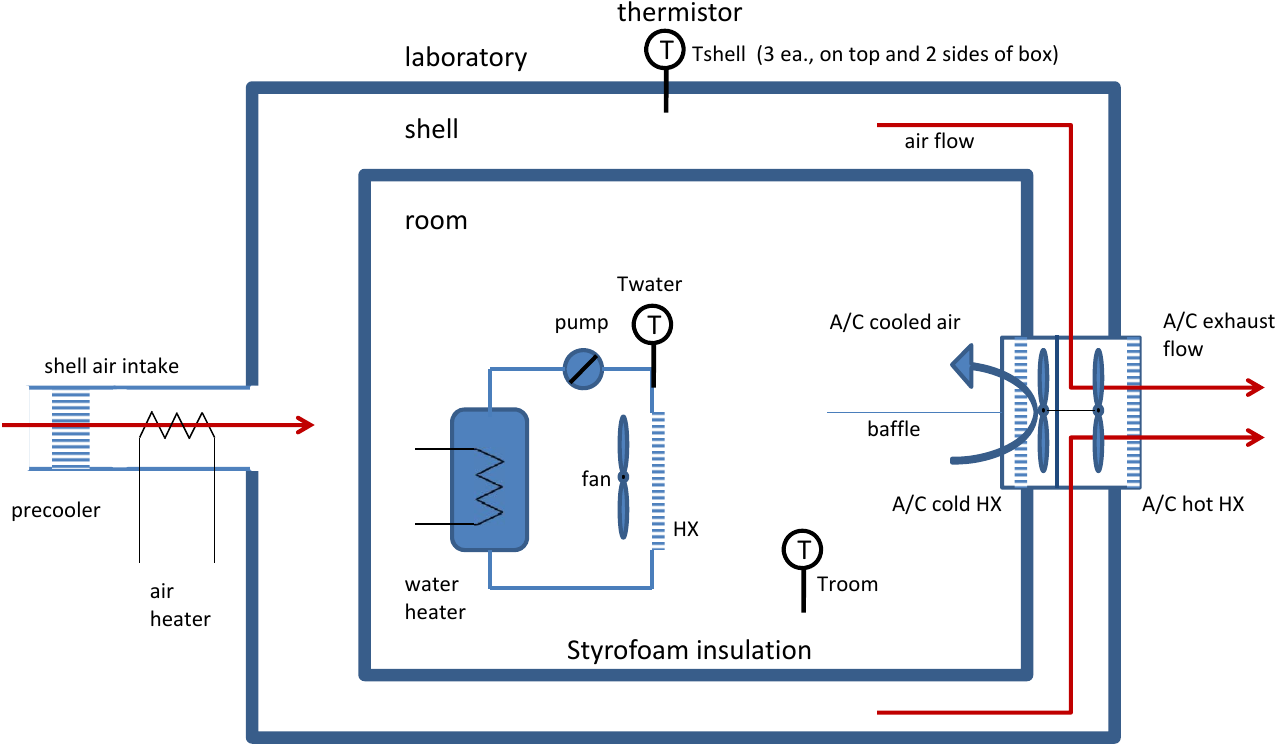}
	\caption{A schematic drawing of the experimental apparatus. Thermistor temperature sensors are indicated by a ``T" inside of a black circle. The insulated walls of the inner enclosure and the components inside the inner enclosure physically emulate the thermodynamics, heat transfer, and air conditioning of a small room.  The water in the water heater emulates the solid heat capacity of the room's walls and floor and other material in the room. The electrical heater in the water heater emulates the heat leak through the room's walls.  Circulating the water through a heat exchanger with fan-driven room air emulates the heat transfer from the room's walls and floor and is tuned to provide a realistic air conditioner (AC) duty cycle. The shell space created by the insulated outer enclosure is fed from the lefthand side with temperature regulated air to shield the inner enclosure from temperature variability in the laboratory. The air flow in the shell is exhausted though the AC's hot heat exchanger and provides constant temperature intake air to isolate the AC operation from the temperature variability in the laboratory. The AC's thermostat is replaced with a software controlled relay that enables modulation of the AC set point temperature or other advanced control functions. }
	\label{fig:apparatus}
\end{figure*}

The experimental test bed for this work is a physically emulated, single-room office as shown in Fig.\ \ref{fig:apparatus}.  The test bed consists of two roughly cubic, nested enclosures. The inner and outer enclosures are constructed from Dow Scoreboard\cite{dow_insulation} insulating foam sheets.  The inner enclosure is roughly 1.3 m on a side and includes a plywood floor to distribute the weight of equipment so as not to damage the foam insulating floor.  Inside the inner enclosure is a 20 gallon, point-of-use water electric heater. The water in the tank comprises most of the heat capacity in the inner enclosure and emulates the heat capacity of the solid materials typically found in an office or room.  Although the water heater is rated to 1000 W, it is driven with a programmable solid state relay (SSR) to provide 200 W on average, which is about half of the total heat dissipated in the inner enclosure. The other heat dissipation is contributed by the electrical loads of a fan and a pump, both described below.  The heat generated in the inner enclosure is in turn removed by a simple window-mount AC that is located with its grille nearly flush with one wall in the inner enclosure.  A cardboard baffle is placed perpendicular to the face of the air conditioner between the cool air supply louvers and the return flow vents to ensure that the cool air flow is effectively mixed with the rest of the air in the inner enclosure and is not "short-circuited" to the return by the small dimensions of the inner enclosure.

In a typical office, there are numerous solid objects, including the drywall bounding the space, that have substantial surface area and heat capacity such that the total heat capacity that the air conditioner works to cool in each cycle is much larger than that just of the air itself. The thermal resistance between the air and these solids creates a second-order heat transfer model \cite{Zhang_2013} where the air cools faster than the surrounding solids and the solid-to-air heat transfer coefficient determines the rate at which the solids cool. Since the water tank is insulated and of compact shape, it was necessary to provide sufficient thermal contact between the water and the air in order to emulate the solid-to-air heat transfer in a typical room. The heat transfer is controlled by using a hydronic loop and pump to circulate the warm water from water tank through a metal-fin heat exchanger.  A fan (separate from two fans in the AC) is mounted on one face of the heat exchanger and vigorously circulates the inner enclosure air through it, providing thermal contact between the heated water and the air.  This fan also stirs the air in the office sufficiently to prevent hot or cold spots from developing due to stagnation.  Initially, the power consumption of the water pump (75 W) and the fan (133 W) were assumed to be dissipated directly to the air. Subsequent modeling and system identification showed that the heat generated in these solid objects can be treated approximately as if all the heat were deposited in the water instead.

To explore the dynamic effects of thermostatic control on the AC, the built-in, bulb-type mechanical thermostat was replaced with instrumentation monitored and controlled by a computer running custom LabView software.  In addition to acting as a smart thermostat, the software extracts and logs temperature measurements from thermistors attached to the apparatus, and it allows the operator to specify several parameters including the temperature set point, the dead band, and independent time-dependent modulation of the set point and dead band edges. The software also enables direct control over the on/off state of the AC, but this capability is not used in this work.  One thermistor monitors the air temperature inside the inner enclosure and is used as the input to the LabView-based thermostatic control.  The software uses this reading to turn on or off the AC's compressor through a relay.  Other thermistors are used to measure the temperature of the laboratory, the temperature of the water in the hydronic loop, and various environmental temperatures surrounding the office.

The AC itself is set to a mode in which its fan runs continuously, which helps mix the air inside the inner enclosure and simultaneously ensures a continuous flow of constant-temperature air across the AC's hot heat exchanger.  The total power consumed by the AC's refrigerant compressor and fans is monitored by a Hameg Model HM 8115-2 autoranging power meter, and the power data are read by the LabView software through a USB interface.

In early experimental trials, we observed that the 2 in.\ thick styrofoam insulation of the inner enclosure was not sufficient to eliminate effects from the diurnal cycle of temperature in the laboratory, which could be as large as 8 $^\circ$C over a single day.  A second outer enclosure, the "shell" of Fig.\ \ref{fig:apparatus}, was built surrounding the experiment in order to control the temperature outside the walls of the inner enclosure as well.  Air from the laboratory is first precooled by a chilled water loop and then passed through an air heater that is controlled by the computer.  The computer monitors temperature on three sides of the shell and modulates the electrical heat following the intake air precooler. This technique maintains the average temperature measured by the three thermistors in the shell space constant to less than 0.3 $^{\circ}$C. The temperature in the shell is chosen to be equal to the air temperature set point in the inner enclosure so that, on average, the heat flow through the inner enclosure walls is minimized.  We also note that the intake air for the AC's hot heat exchanger is drawn from the temperature regulated shell space. This configuration stabilizes the ``external environmental" conditions for the AC operation. By stabilizing the conditions for heat flow into the inner enclosure and for AC operation, the time-dependent power consumption of the AC can be studied without the environmental variabilities that plague such measurements on real buildings.

\begin{figure}[thpb]
	\centering
	\includegraphics[width=0.4\textwidth]{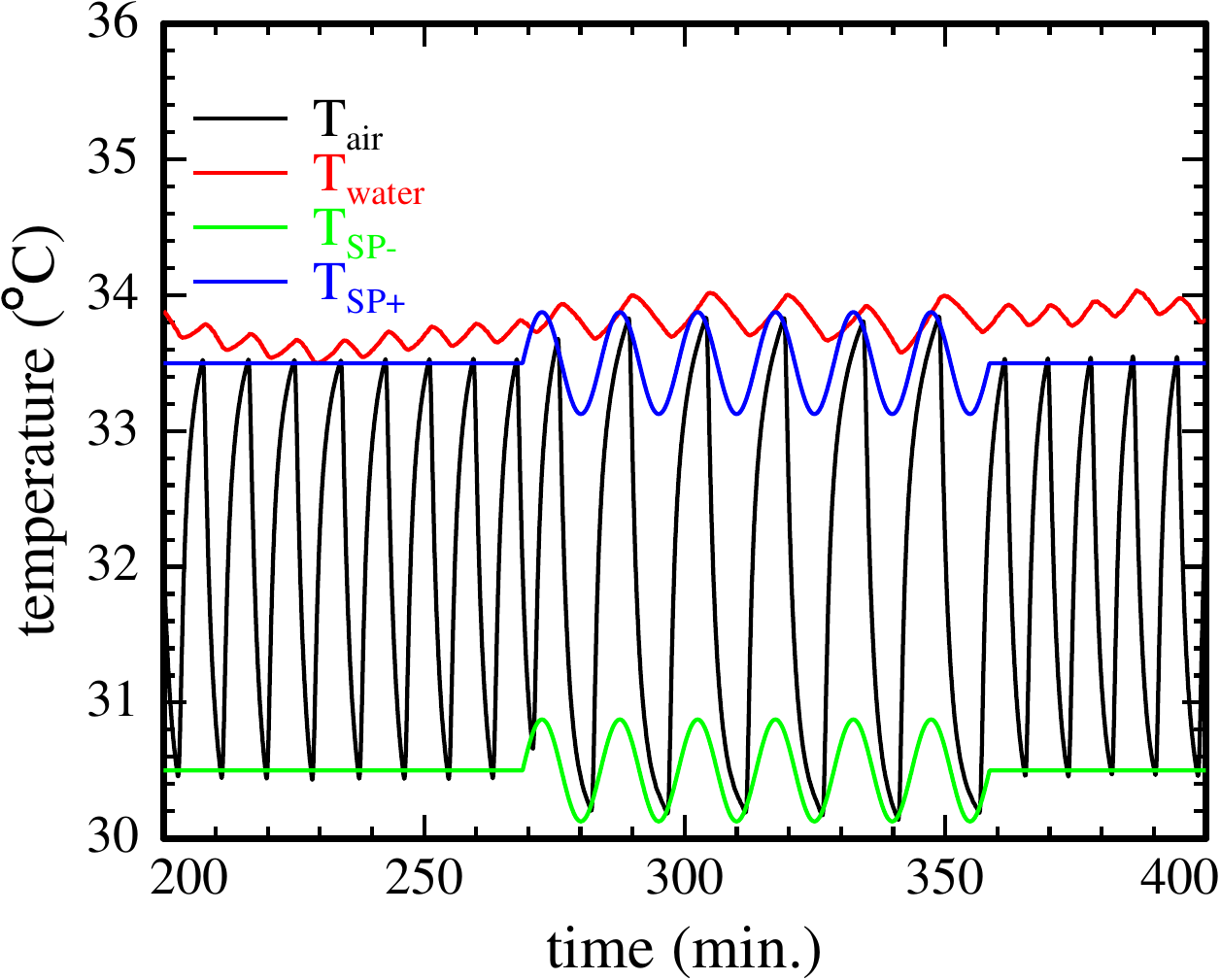}
	\caption{A sample of the raw temperature data for the inner enclosure air temperature $T_{\rm air}$, the circulating water temperature $T_{\rm water}$, and the upper and lower dead band limits, $T_{\rm SP+}$ and $T_{\rm SP-}$, respectively. The temperature set point (and dead band limits) is modulated sinusoidally for an interval of about 90 minutes. The modulation period and amplitude are parameters defining the experiment. The modulations are cycled on and off every 90 minutes to allow for frequent reassessment of the AC baseline operation.}
	\label{fig:rawdata}
\end{figure}

Even with these efforts to create a well stabilized thermodynamic control volume around the inner enclosure, there are still some sources of systematic error from the environment and from the off-the-shelf equipment itself.  For example, there is a residual effect of the ambient temperature of the laboratory on the AC behavior (both the duty cycle and the average power consumption), which may be due to the finite flow rate of the temperature regulated air in the shell through the hot heat exchanger of the AC.  In order to reduce these residual nuisance effects over the multiple weeks of measurements required to fully map the control response functions and cost of control, the data acquisition software alternates intervals of set point modulation with intervals of fixed set point. This experimental procedure continually re-establishes the baseline power consumption of the AC and provides a differential measurement of AC natural period and power consumption.  An example chart of the temperature data for such an experiment is shown in Fig.\ \ref{fig:rawdata}.

With a constant heating rate into the water, the experimental system locks to a specific phasing with respect to the imposed temperature modulations, in the absence of any other perturbations to the inner enclosure.  This is an unrealistic situation and may lead to locally stable states that are highly dependent on the particular heat transfer characteristics of the experiment and the specific temperature trajectory of the experiment at startup.  In addition, in physically relevant residential or commercial spaces the heat load in the space is randomized by the time variation of occupancy, by the time variation of the activities inside the space, and by the environment.  Therefore, to find the most robust state of the system, the heat input to the water was randomized by adding zero-mean noise to the time-averaged heating rate of the water heater.  For the experiments shown, this noise consists of random step changes between $\pm100$ W to the water heater, applied for random intervals averaging to about 15 minutes.  The software allows control of the heating rate, the deviation from the mean rate, and the reset time; and it controls the physical heater via an output voltage to a circuit that converts the signal to 4-20 mA for programming the SSR.  The SSR itself passes 120 V AC power, stepped down by a variac, to the heater element.
\section{Results}\label{sec:results}

Using the test bed described above, we systematically map out the electrical power response of the window-mount AC to sinusoidal perturbations of the set point and dead band edges of the form
$$
T_{\rm set point}\left( t \right) = T_{\rm avg} +\Delta T \sin \left( 2 \pi t/ \tau \right). \eqno{(1)}
$$
Experiments are performed with temperature modulation amplitudes from $\Delta T$=0.0625 to 0.5 $^{\circ}$C and modulation periods ranging from $\tau$=3 to 40 minutes per cycle.  Although not all experiments are the same length, most experiments were run for about 3.5 days---long enough to ensure that any transients at the beginning of the experiment have settled (1-3 hours for starting a new experiment) and do not significantly affect the overall results. This duration provides a large number of modulation cycles, 62 even for the longest $\tau$=40-minute modulation period.

\begin{figure*}[thpb]
	\centering
	\includegraphics[width=0.95\textwidth]{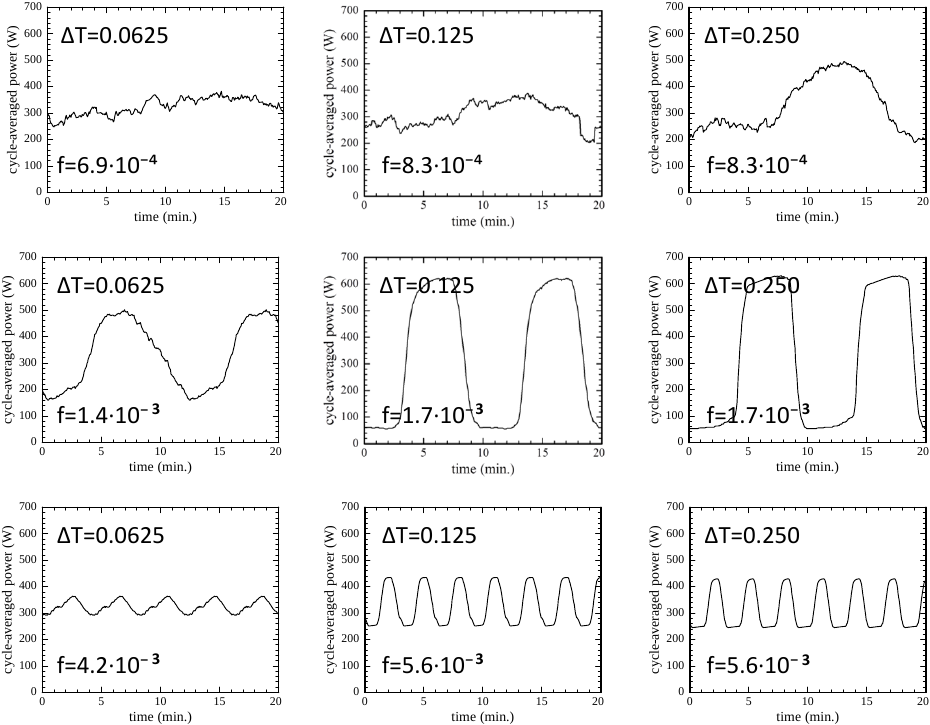}
	\caption{Cycle-averaged AC electrical power during modulation of the temperature set point by three different amplitudes, $\Delta T$=0.0625, 0.125, and 0.25 $^\circ C$, and three different periods, $\tau$=3, 10 and 20 minutes.  The periods in the first column ($\Delta T$=0.0625 $^{\circ}C$) are slightly different from those in the next two columns due to changes in the experiment implemented between these data sets.}
	\label{fig:respfuncplots}
\end{figure*}

For each experiment at a given $\Delta T$ and $\tau$, the beginning of each modulation cycle is identified and the AC's electrical power at each time point in the cycles is averaged to develop a single, average electrical power response to the sinusoidal set point perturbation. The AC's electrical power response for several experimental conditions are shown in Fig.~\ref{fig:respfuncplots}. The columns in Fig.~\ref{fig:respfuncplots}   are each for a single value of amplitude, $\Delta T=$ 0.0625, 0.125, and 0.25 $^{\circ}C$ from left to right.  The plots in each row have nearly the same modulation period, $\tau =$ 20, 10, and 3 minutes from top to bottom.  (The actual values for the leftmost column are 24, 12, and 4 minutes.  These data were taken in an earlier set of experiments on the same apparatus.)

The responses to the harmonic modulations are not simple sinusoids, even for the lowest amplitude in Fig.~\ref{fig:respfuncplots}, as the finite temperature dead band creates significant nonlinearity.  The nonlinearity becomes more pronounced for larger set point modulation amplitudes.  Near the natural period of the unperturbed system (about 9 minutes), the response function at high modulation amplitudes almost resembles a square wave. The square wave is a result of saturation because the AC cannot provide greater peak-to-peak changes in power than the maximum power it consumes. At longer modulation periods, the nonlinearity also generates frequency doubling with the AC's electrical power response having significant variability at half the modulation period.

The responses in Fig.~\ref{fig:respfuncplots} are in contrast to those from an HVAC in a large commercial building \cite{beil}. In those systems, damper levels and fan speeds are continuously adjusting to regulate around the temperature set point. The simpler controls on a window-mount AC unit can only turn the compressor on or off and cannot regulate temperature at all within the dead band.

\begin{figure}[thpb]
	\centering
	\includegraphics[width=0.4\textwidth]{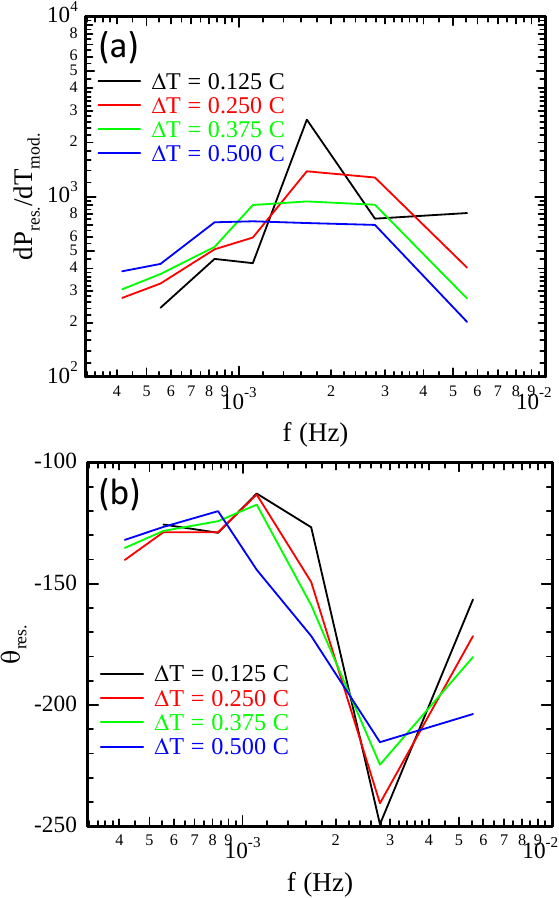}
	\caption{Frequency response of the AC electrical power for set point modulation amplitudes of $\Delta T$=0.125, 0.250, 0.375, and 0.500 $^\circ C$ :  (a) Amplitude response (b)  Phase shift.}
	\label{fig:bodeplots}
\end{figure}

For each plot Fig.~\ref{fig:respfuncplots}, the magnitude of the variation in consumed power is a measure of the amount of power that the individual ACs can provide for regulation services with a given timescale and an allowed temperature excursion in the room. Although the response is not purely sinusoidal, it is cyclic and can be characterized by extracting the magnitude and phase of the first harmonic of cycle-averaged electrical power waveform.  These modulation amplitude-dependent response functions are shown in Fig.~\ref{fig:bodeplots} as functions of modulation frequency.

Modulation of the set point of the window-mount AC unit does not only modulate the electrical power, but it also creates an increase in the average electrical power consumption reflecting the cost of control. The increase in consumption is difficult to measure because, at least for the lower temperature modulations studied ($\Delta T \le 0.125$ C), the loss is comparable to other sources of noise in the system.  By alternating periods of modulation with quiescent periods, as shown in Fig.\ \ref{fig:rawdata}, it is possible to reduce these systematic errors and obtain values reproducible to less than 5 W.  The extra AC electrical power consumption $\Delta P$ is calculated from
$$
\Delta P = (P_{\rm mod} - \dot{Q}_{\rm w, mod}/\eta) - (P_{0} - \dot{Q}_{\rm w, 0} /\eta) \eqno{(2)}
$$
where $P_{\rm mod}$ and $P_0$ are the time-averaged AC power with and without modulation, respectively; $\dot{Q}_{\rm w, mod}$ and $\dot{Q}_{\rm w, 0}$ are the time-averaged electrical heat applied to the water heater with and without modulation, respectively; and $\eta$ is the coefficient of performance of the air conditioner.

\begin{figure}[thpb]
	\centering
	\includegraphics[width=0.4\textwidth]{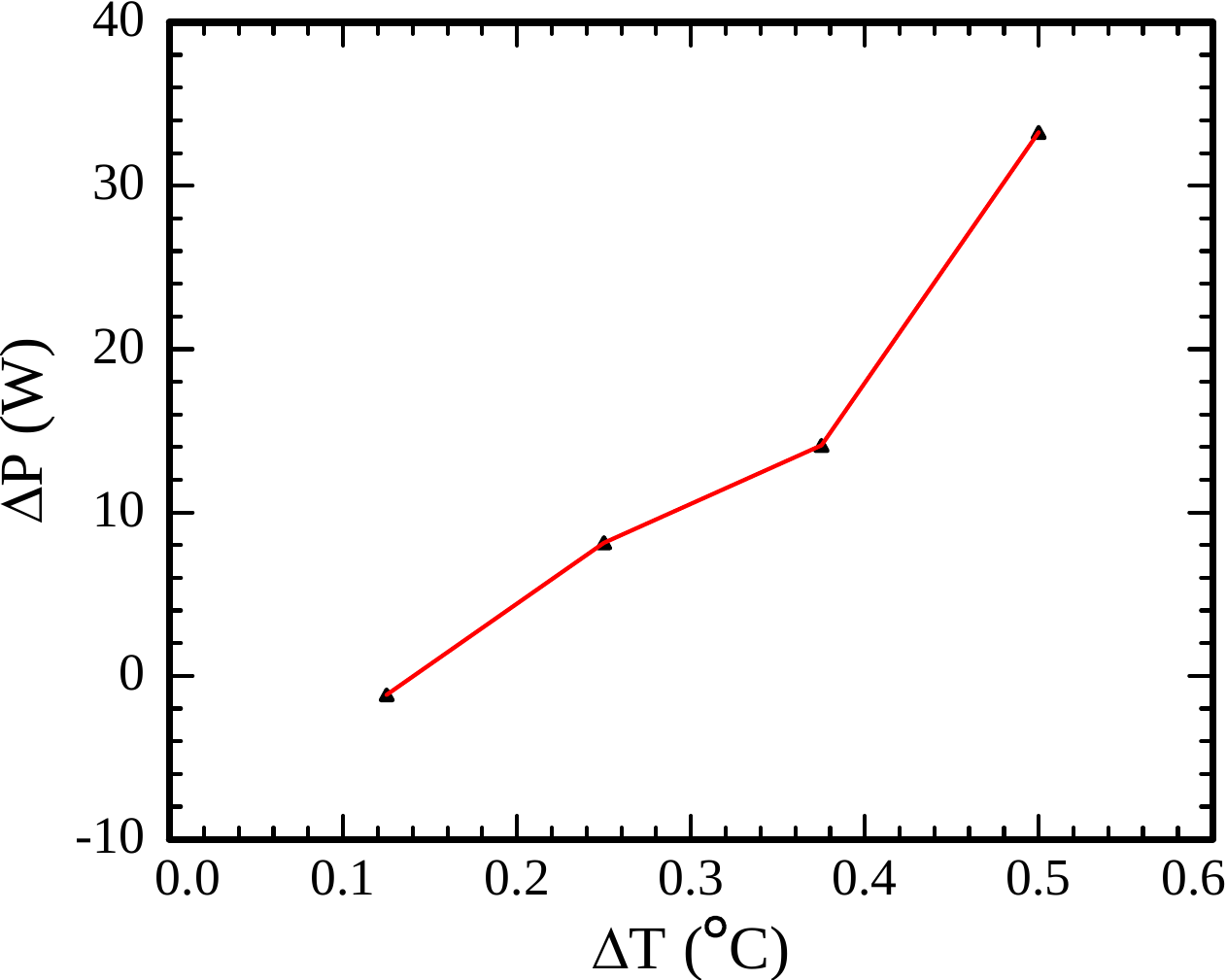}
	\caption{The excess power consumption $\Delta P$ of the air conditioner attributable to temperature set point modulation of amplitude $\Delta T$.  The period of the modulation is $\tau$=20 minutes.}
	\label{fig:dissipation}
\end{figure}

The extra AC electrical power consumption is shown in Fig.\ \ref{fig:dissipation} as a function of temperature amplitude for a modulation period of 20 minutes. In Fig.~\ref{fig:dissipation}, $\Delta P$ extrapolates to a value less than zero, which cannot be consistent with definition of $\Delta P$. This inconsistency may be due to unresolved systematic errors in our experimental methods. Alternatively, $\Delta P$ may actually be zero for $\Delta T <$0.125 $^\circ C$.

\section{Discussion}\label{sec:discussion}

The results up to this point have been presented in terms of the response of a single window-mount AC subject to set point modulation. However, the thermal noise discussed in Section~\ref{sec:testbed} allows our results to be interpreted as the response of a large ensemble of window-mount ACs.  The thermal noise applied to the circulating water inhibits the dynamical response of the single AC from locking onto the perturbation frequency and forces the single AC to more fully explore the phase space of temperatures and on/off states. Over many cycles for a modulation amplitude and period pair, the single AC starts each cycle in a different state. We can then interpret cycle averaging of the individual modulation cycles as being equivalent to ensemble averaging over a large number of ACs that start each modulation cycle at a random point in their duty cycle. Interpreted in this way, the results in Figs.~\ref{fig:respfuncplots}-\ref{fig:dissipation} are more broadly applicable to large ensembles of window-mount AC units.

The shape of the frequency response in Fig.~\ref{fig:bodeplots} is in many ways similar to the $H_1$ response in Fig. 8 of \cite{lin}, but the physical causes of the response are often different. Considering the lower amplitude modulations of $\Delta T$ = 0.125 or 0.250 $^\circ C$, the results in Fig.~\ref{fig:bodeplots} show a similar peak at a mid range of frequency with a roll off at both higher and lower frequencies. The peak at approximately 1 minute and high frequency rolloff in \cite{lin} is a result of the dynamics of the AHU fan, which may be natural dynamics or acceleration limiters in the variable frequency drive that ultimately controls the fan speed. This peak is not expected to move as the conditions in the building or the environmental conditions change.

In Fig.~\ref{fig:bodeplots}, the peak in the response magnitude at approximately 9 minutes is associated with the natural duty cycle of the AC unit. If the set point and dead band limits are modulated much faster than the natural duty cycle, the set point explores the high and low extremes of its modulation many times during an AC duty cycle. The AC unit will always switch on/off states at these ``inner'' extremes of the dead band modulation rendering the control input ineffective in modulating the average power consumption. This effect results in the roll off of the response at high frequency in Fig.~\ref{fig:bodeplots}. We note that the natural duty cycle and the location of the peak (and $\approx$ 180$^\circ$ phase shift) in Fig.~\ref{fig:bodeplots} is expected to move as conditions change creating challenges for designing a closed loop controller for an ensemble of window-mount ACs.

The low-frequency roll off in the response in Fig.~\ref{fig:bodeplots} is caused by the natural adjustment of the average air and solid (i.e. water) temperature to slow modulations in set point. Under these slow modulation conditions, the average temperature quickly adjusts to the new set point and the effective ``error'' between the set point and average temperature is reduced. Without significant temperature error, the AC's electrical consumption tends toward its baseline consumption and the control signal does not provide significant modulation of AC electrical power. Related slow adjustments of the system state in \cite{lin} occur to reduce the level of electrical power modulation.

The experiments in \cite{lin} purposefully limited the range of the control inputs to avoid possible damage to the fan and connected ductwork, and therefore system identification may not have encountered significant nonlinearities. In contrast, the present work encounters nonlinearities even at relatively small set point modulations. The gradual reduction in the peak of the magnitude response at the fundamental frequency in Fig.~\ref{fig:bodeplots} is caused by the saturation in the response of the AC unit. This can be seen in the middle row of Fig.~\ref{fig:respfuncplots}. There is an associated impact on the phase shift at the fundamental frequency in Fig.~\ref{fig:bodeplots}. The presence of significant nonlinearity at even $\Delta T \approx$0.250 $^\circ C$ is important because the set point temperature resolution for most thermostats is at best 0.1$^\circ C$.

The frequency response functions in Fig.~\ref{fig:bodeplots} may be used to design a control that accurately tracks a frequency regulation reference signal.  However, the extra time-average power consumption in Fig.~\ref{fig:dissipation} is key to understanding if providing ancillary services such as frequency regulation is economically viable.  In Fig.~\ref{fig:dissipation}, a modulation of $\Delta T$=0.250 $^\circ C$ at a $\tau$=20 minute period would cause about 10 Watts of extra energy consumption. One could estimate from Figure~\ref{fig:respfuncplots} that this same set point modulation would generate roughly 100 Watts of zero-to-peak response in the cycle-averaged power.  Assuming that the AC owner pays retail for energy at \$200/MW-hr and is paid an optimistic frequency regulation market clearing price of \$40/(MW/hr), the ratio of cost of energy to frequency regulation income is 1:2. Fifty percent of the AC owners' ancillary service revenue is lost to extra energy consumption. The relatively close balance of income and cost of providing the ancillary service clearly points to the need for more extensive measurements of $\Delta P$ under ideal settings, as performed here, and under realistic settings using historical frequency regulation signals.

\section{Conclusion}\label{sec:conclusions}

In this manuscript, we have demonstrated experimental methods to characterize the frequency response and cost of control of window-mount air conditioning units. These methods are extensible to characterization of many other forms of thermostatically controlled loads (TCL).  By injecting noise into the system under test, we have made a preliminary demonstration of how to measure the frequency response of an ensemble of TCLs using a single TCL. The value of these methods is twofold. First, experimental response functions capture crucial unmodeled behavior of loads and other distributed energy resources before deploying control systems to manage them in real power grids. Second, they enable the exploration of response functions of large ensembles without having to incur the expense of building or deploying large numbers of devices.

There are many directions for future work using the same methods, including:
\begin{itemize}
\item Further investigation of the noise injection method used in this work to understand the accuracy of the recovered response function for ensembles of devices 

\item Extending the frequency domain of the window unit air conditioner characterization to better capture the control response for faster time scale ancillary services like RegD in the PJM market\cite{PJM1}

\item Improving and expanding the excess energy measurements to better characterize the economic viability of providing frequency regulation from an ensemble of small thermostatically controlled loads 

\item Devising related measurement schemes to isolate and understand loss mechanisms in thermostatically controlled loads to better identify control methods that minimize the losses

\item Extending the measurements performed here to include characterization of the state probability distribution functions, a characterization potentially more useful to more advanced control algorithms
\end{itemize}

Finally, it is not possible to test control algorithms for ensembles of thermostatically controlled loads using these methods on a single TCL. The time averaging used to capture the ensemble behavior does not allow for observation of the ensemble response in real-time. To test any control algorithm developed using these experimental methods requires significantly scaling up the number of these test beds to adequately represent ensemble behavior in real time.

\addtolength{\textheight}{-12cm}   



%

\section*{ACKNOWLEDGMENT}
This work was supported by the Microgrid Program of the Office of Electricity within the U.S. Department of Energy. We are thankful to Gregory Swift for fruitful technical discussions and to Carmen Espinoza for his skill in constructing the experimental test bed.  Research was carried out under the auspices of the National Nuclear Security Administration of the U.S. Department of Energy at Los Alamos National Laboratory under Contract No. DE C52-06NA25396. 



\bibliographystyle{unsrt}
\bibliography{HVAC}

\end{document}